\documentclass{article}
\usepackage{amsthm}

\theoremstyle{plain}
\newtheorem{theorem}{Theorem}[section]
\newtheorem{proposition}[theorem]{Proposition}

\newtheorem{example}[theorem]{Example}
\newtheorem{definition}[theorem]{Definition}

\usepackage{multirow} 
\usepackage{enumitem} 
\usepackage{xspace} 

\usepackage{amsmath,amssymb}
\usepackage{bm}     
\usepackage{pdfsync}

\usepackage[T1]{fontenc}

\def\Mat{\mathop{\rm Mat}\nolimits}

\def\gcd{\mathop{\rm gcd}\nolimits}

\newcommand \ie {\textit{i.e.}}
\newcommand \eg {\textit{e.g.}}

\DeclareMathOperator\IntrinsicContent{FD}
\newcommand \prim{\mathop{\rm prim}\nolimits}
\newcommand \FF {{\mathbb F}}
\newcommand \NN {{\mathbb N}}
\newcommand \QQ {{\mathbb Q}}

\newcommand \ZZ {{\mathbb Z}}
\newcommand \CC {{\mathbb C}}

\newcommand \TTo[1]{\mathop{\longrightarrow}\limits^{#1}}

\newcommand \longiso{\,\smash{\TTo{\lower 7pt\hbox{$\scriptstyle\sim$}}}\,}

\newcommand\MT{M\"obius transformation\xspace}
\newcommand\MTs{M\"obius transformations\xspace}

\begin{document}

\title{Certifying Irreducibility in {$\ZZ[x]$}}

\author{John Abbott}


\maketitle

\begin{abstract}
We consider the question of certifying that a polynomial in $\ZZ[x]$
or $\QQ[x]$ is irreducible.  Knowing that a polynomial is irreducible
lets us recognise that a quotient ring is actually a field extension
(equiv.~that a polynomial ideal is maximal).   Checking that a polynomial
is irreducible by factorizing it is unsatisfactory because it requires
trusting a relatively large and complicated program (whose correctness
cannot easily be verified).  We present a practical method for generating
certificates of irreducibility which can be verified by relatively
simple computations; we assume that primes and irreducibles in $\FF_p[x]$
are self-certifying.
\end{abstract}




\section{Introduction}
\label{sec:Intro}

\subsection{What is a ``certificate''?}

A certificate that object $X$ has property $P$ is a ``small'' amount
of extra information $C$ such that some \textit{quick and simple computations}
with $X$ and $C$ suffice to confirm that $X$ does have the property.
We illustrate this vague definition with a well-known, concrete example.

\begin{example}
  We can certify that a positive integer $n$ is prime using a
  \textit{Lucas-Pratt certificate}~\cite{zbMATH03495595}.  The idea is to find a witness $w$
  such that $w^{n-1} \equiv 1 \mod n$ and $w^{(n-1)/q} \not\equiv 1 \mod
  n$ for all prime factors $q$ of $n{-}1$.
  
  These certificates have a recursive structure, since in general
  we must certify each prime factor $q$ of $n{-}1$.  To avoid
  infinite recursion we say that all small primes up to some limit
  are ``self-certifying'' (\ie~they need no certificate).

  Thus a Lucas-Pratt certificate comprises a witness $w$, and a list of
  prime factors $q_1,q_2,\ldots$ of $n{-}1$ (and certificates for each $q_j$).
  Verification involves:
  \begin{itemize}
    \setlength{\itemsep}{1pt}
  \item verify that $w^{n-1} \equiv 1 \mod n$;
  \item verify that each $w^{(n-1)/q_j} \not\equiv 1 \mod n$;
  \item verify that $n{-}1 \,=\, \prod_j q_j^{e_j}$ for positive exponents $e_j$;
  \item recursively verify that each $q_j$ is prime.
  \end{itemize}

  The operations required to verify such a certificate are:
  iteration over a list, exponentiation modulo an integer,
  comparison with~$1$, division of integers, and divisibility testing
  of integers.  These are all simple operations, and the entire
  function to verify a Lucas-Pratt certificate is small enough to be
  fully verifiable itself.
\end{example}

An important point in this example is that the certificate actually
involves several cases: namely, if the prime is small enough, the certificate
just says that it is a ``small prime'' (\eg~we can verify by table-lookup); otherwise the certificate contains
a non-trivial body.  In this instance there are just two possible cases.

We note that \textit{generating} a Lucas-Pratt certificate could be costly because
the prime factorization of $n{-}1$ must be computed.

\subsection{Costs of a certificate}

The total cost of a certificate comprises several components:
\begin{itemize}
\setlength{\itemsep}{0pt}
\item computational cost of generating the certificate;
\item size of the certificate (\eg~cost of storage or transmission);
\item computational cost of verification given the certificate;
\item size and code complexity of the verifier.
\end{itemize}

In the case of certifying the irreducibility of a polynomial in $\ZZ[x]$
we could issue trivial certificates for all polynomials, and say that
the verifier simply has to be an implementation of a polynomial factorizer.
We regard this as unsatisfactory because the size and code complexity of
the verifier are too high.

\section{Irreducibility Criteria for $\ZZ[x]$ and $\QQ[x]$}

We can immediately reduce from $\QQ[x]$ to $\ZZ[x]$ thanks to Gauss's
Lemma (for polynomials): let $f \in \QQ[x]$ be non-constant then $f$
is irreducible if and only if $\prim(f) \in \ZZ[x]$ is irreducible,
where $\prim(f) = \alpha f$ and the uniquely defined, non-zero factor
$\alpha \in \QQ$ is such that all coefficients of $\prim(f)$ are integers
with common factor~$1$, and the leading coefficient is positive,

\medskip

The problem of certifying irreducibility in $\ZZ[x]$ has a long history,
and has already been considered by several people.
Here is a list of some approaches:
\begin{itemize}
\setlength{\itemsep}{0pt}
\item give a ``large'' evaluation point $n$ such that $f(n)$ has a large prime factor;
\item degree analysis (from factorizations over one or more finite fields)\footnote{\textit{degree analysis} has likely been known for a long time};
\item a linear polynomial is obviously irreducible;
\item Newton polygon methods (\eg~Sch\"onemann, Eisenstein, and Dumas~\cite{zbMATH02646458});
\item Vahlen-Capelli lemma~\cite{zbMATH03204669} for binomials 
\item Perron's Criterion~\cite{zbMATH02643268}; 
\item the coefficients are (non-negative) digits of a prime to some base $b$ (\eg~\cite{zbMATH02643268}).
\end{itemize}
The first technique in the list was inspired by ideas from~\cite{zbMATH03974286}; it seems to be new.

In this presentation, we shall assume that the degree is at least $2$, and
shall concentrate on the first two methods as they are far more widely
applicable than others listed.

\subsection{Factor Degree Analysis}

Factor degree analysis is a well-known, behind-the-scenes technique in
polynomial factorization.  It involves using degrees of modular
factors to obtain a list of \textit{excluded degrees} for factors in
$\ZZ[x]$.

We define a \textbf{factor degree lower bound} for $f \in \ZZ[x]$ to
be $\Delta \in \NN$ such that we have excluded all degrees less than
$\Delta$, \eg~through factor degree analysis.  We can certify this
lower bound by accompanying it with the modular factorizations used.
Clearly, if degree analysis excludes all degrees up to $\frac{1}{2}
\deg f$ then \textit{we have proved that $f$ is irreducible.}  
Finally, we may always take $\Delta=1$ without any degree analysis.


In many cases we can indeed prove/certify irreducibility via degree analysis.
However, there are some (infinite) families of polynomials where one
must use ``larger'' primes, and there are also (infinite) families
where irreducibility cannot be proved via factor degree analysis
(\eg~resultants, in particular Swinnerton-Dyer polynomials, see
also~\cite{zbMATH03841227}).

\begin{example}
  The well-known, classical example of a polynomial which cannot be
  proved irreducible by degree analysis is $x^4+1$: every modular
  factorization is into either $4$ linears or $2$ quadratics, so this
  does not let us exclude the possible existence of a degree~$2$
  factor.

  There are also many polynomials which can be proved irreducible by
  degree analysis, but are not irreducible modulo any prime; this
  property depends on the Galois group of the polynomial.  For
  instance, $f = x^4+x^3+3x+4$ is one such polynomial: modulo~2 the
  irreducible factors have degrees~1 and~3, and modulo~5 both factors
  have degree~$2$; but it is never irreducible modulo $p$.
\end{example}

\subsubsection{Degree analysis certificate}

A degree analysis certificate comprises
\begin{itemize}
  \item a subset $D \subseteq \{1,2,\ldots,\frac{1}{2}\deg f\}$ of ``not excluded'' factor degrees
\item  a list, $L$, of pairs: a prime $p$, and the irreducible factors of $f$ modulo $p$
\end{itemize}
If $D = \emptyset$, we have a certificate of ireducibility;
otherwise the smallest element of the set is a \textit{factor degree lower bound.}

Verification of the certificate involves the following steps:
\begin{itemize}
  \setlength{\itemsep}{0pt}
\item for each entry in $L$, check that the product of the modular factors is $f$;
\item for each entry in $L$, compute the set of degrees of all possible products of the modular factors; verify that their intersection is $D$;
  \item check that each modular factor is irreducible (\eg~use gaussian reduction to compute the rank of $B{-}I$ where~$B$ is the Berlekamp matrix).
\end{itemize}
The main cost of the verification is the computation of $B$ and the
rank of $B{-}I$; the cost of computing $B$ is greater for larger
primes, so we prefer to generate certificates which use smaller primes
if possible.

\subsubsection{Practical matters}

We would like to know, in practice, how costly it is to produce a
useful degree analysis certificate, and how large the resulting
certificate could be.  More specifically:
\begin{itemize}
  \setlength{\itemsep}{0pt}
\item How many different primes should we consider?  And how large?
\item How to find a minimal set of primes yielding the factor degree subset?
\item How many primes are typically in the minimal set?
\end{itemize}
In our experience, a minimal length list very rarely contains more than~$3$
entries, but we should expect to consider many more primes during
generation of the certificate.  We can construct irreducible
polynomials which require considering ``large'' primes to obtain
useful degree information (\eg~$x^2+Nx+N$ where $N = 1000000!$) but in
many cases ``small'' primes up to around $\deg f$ suffice.

\subsection{Irreducibility Certificates for $\ZZ[x]$ via Evaluation}
\label{subsec:IrredByEvaln}

Bunyakowski's conjecture (\eg~see page~323 of~\cite{zbMATH00233957}) states that if $f \in \ZZ[x]$ is irreducible
(and has trivial fixed divisor) then $|f(n)|$ is prime for infinitely
many $n \in \ZZ$.  Assuming the conjecture is true, we can get a
certificate of irreducibility by finding a suitable evaluation point
$n$ (and perhaps including a certificate that $|f(n)|$ is prime).

Applying Bunyakowski's conjecture directly is inconvenient for two reasons:
\begin{itemize}
\setlength{\itemsep}{1pt}
\item we want to handle polynomials with non-trivial fixed divisor;
\item finding a suitable $n$ may be costly, and the resulting $|f(n)|$ may be large.
\end{itemize}
The first point is solved by an easy generalization of the conjecture:
let $f \in \ZZ[x]$ be irreducible and $\delta$ be its fixed divisor,
then there are infinitely many $n \in \ZZ$ such that $|f(n)|/\delta$ is prime.
The second point is a genuine inconvenience: for some polynomials, it can
be costly to find a ``Bunyakowski prime,'' and the prime
itself will be large (and thus costly to verify).  For example, let $f = x^{16}+4x^{14}+6x^2+4$
then the smallest good evaluation point is $n=6615$, and $|f(n)| \approx 1.3 \times 10^{61}$.

\subsubsection{A large prime factor suffices}
\label{subsubsec:PrimeFactorSuffices}

Here we present a much more practical way of certifying irreducibility by evaluation:
we require just a sufficiently large prime factor.
Let $f \in \ZZ[x]$ be non-constant, and let $\rho \in \QQ$ be a \textbf{root
bound} for $f$: that is, for every $\alpha \in \CC$ such that
$f(\alpha)=0$ we have $|\alpha| \le \rho$.  We note that it is
relatively easy to compute root bounds (\eg~see~\cite{zbMATH06124246}).
The following proposition was partly inspired by Theorem~2 in~\cite{zbMATH03974286}, but appears to be new.

\begin{proposition}
\label{prop:LPFW}
Let $f \in \ZZ[x]$ be non-constant, and let $\rho \in \QQ$ be a root
bound for $f$.  Let $\Delta \in \NN$ be a factor degree lower bound for $f$.
If we have $n \in \ZZ$ with $|n| > 1+\rho$ such that $|f(n)| = sp$
where $s < (|n|-\rho)^\Delta$ and $p$ is prime then $f$ is irreducible.
\end{proposition}

\begin{proof}
  For a contradiction, suppose that $f=gh \in \ZZ[x]$ is a non-trivial
  factorization.  We may assume that $\Delta \le \deg g \le \deg h$.
  We have $f = C_f \prod_{j=1}^d (x - \alpha_j)$ where $d = \deg f$,
  $C_f \in \ZZ$ is the leading coefficient, and the $\alpha_j$ are the
  roots of $f$ in $\CC$.  We may assume that the $\alpha_j$ are indexed
  so that the roots of $g$ are $\alpha_1,\ldots,\alpha_{d_g}$ where $d_g = \deg g$.

  By evaluation we have $f(n) = g(n) \, h(n)$ with all values in $\ZZ$.
  Also $f(n) \neq 0$ since $|n| > \rho$.  We now estimate $|g(n)|$:
  \[
g(n) \; = \; C_g\, \prod_{j=1}^{d_g} (n- \alpha_j)
  \]
  where $C_g \in \ZZ$ is the leading coefficient.  Each factor in the product
  has magnitude greater than 1, so $|g(n)| \,\ge\, (|n|-\rho)^\Delta \,>\, s$.  Similarly,
  $|h(n)| > s$.  This contradicts the given factorization $f(n) = sp$.
\end{proof}

When we have an evaluation point to which Prop.~\ref{prop:LPFW} applies
we call it a \textbf{large prime factor witness} (abbr.~\textbf{LPFW}) for $f, \rho$ and $\Delta$.
We conjecture that every irreducible polynomial
has infinitely many LPFWs; note that Bunyakowski's conjecture implies this.

\begin{example}
  This example shows that it can be beneficial to look for large prime factor
  witnesses rather than Bunyakowski prime witnesses.

  Let $f = x^{12} +12x^4 +92$ and take $\Delta=1$.  We compute
  $\rho=\frac{7}{4}$ as root bound, and then we obtain a LPFW at $n=5$
  with prime factor $p=81382739$.  In contrast, the smallest
  Bunyakowski prime is $\approx 3.06 \times 10^{41}$ at $n=2865$.
\end{example}

In the light of this example we exclude consideration of
a certificate based on Bunyakowski's conjecture, and consider only
LPFWs.

We prefer to issue an LPFW certificate where
the prime $p$ is as small as ``reasonably possible''.  Our
implementation searches for suitable~$n$ in an incremental way, since
smaller values of $|n|$ produce smaller values of $|f(n)|$, and we
expect smaller values of $|f(n)|$ to be more likely to lead to an
``$sp$'' factorization with small prime factor $p$~---~this is only a
heuristic, and does not guarantee to find the smallest such $p$.  We
look for the factorization $|f(n)| = sp$ by trial division by the
first few small primes (and GMP's probabilistic prime test for $p$).

\subsubsection{LPFW certificate}

An LPFW certificate comprises the following information:
\vspace{-2pt}
\begin{itemize}
  \setlength{\itemsep}{0pt}
\item a root bound $\rho$,
\item a factor degree lower bound $\Delta$ \quad$\longleftarrow$ {\small with degree analysis certificate,}
\item the evaluation point~$n > 1+\rho$,
\item the large prime factor~$p$ of $|f(n)|$ \quad$\longleftarrow$ {\small (opt.) with certificate of primality.}
\end{itemize}
\goodbreak
\noindent
Verification of an LPFW certificate entails:
\vspace{-2pt}
\begin{itemize}
  \setlength{\itemsep}{0pt}
\item evaluating $f(n)$ and verifying that $p$ is a factor;
\item verifying that the discarded factor $s=|f(n)|/p$ satisfies $s  < (|n|-\rho)^\Delta$;
\item verifying that $\rho$ is a root bound for $f$ \quad$\longleftarrow$ see comment below;
\item (if $\Delta > 1$) verifying that $\Delta$ is a factor degree lower bound;
\item verifying that $p$ is (probably) prime.
\end{itemize}
\goodbreak
In many cases the root
bound can be verified simply by evaluation of a modified polynomial: let $f(x) = \sum_{j=0}^d a_j x^j$ and set
$f^*(x) = |a_d|x^d - \sum_{j=0}^{d-1} |a_j|x^j$, then if $f^*(\rho) > 0$ then
$\rho$ is a root bound for $f$.  Some tighter root bounds may require applying an
(iterated) Gr\"affe transform to $f$ first (\eg~see~\cite{zbMATH06124246}).

\begin{example}
This example shows how degree information can be useful in finding a
small LPFW.  Let $f = x^4 -1036x^2 +7744$.  We find that
$\rho=33$ is a root bound.  Without degree information (\ie~taking
$\Delta=1$) we obtain the first LPFW at $n=65$ with
corresponding prime $p=13481269$.  In contrast, from the factorization
of $f$ modulo $3$ we can certify that $\Delta=2$ is a factor degree
lower bound for $f$.  This information lets us obtain an LPFW
at $n=47$ with far smaller corresponding prime $p=14519$.
\end{example}

\section{M\"obius Transformations}

We define a (minor generalization of) a \textit{\MT} for $\ZZ[x]$.
The crucial property for us is that these transformations preserve
irreducibility (except for some polynomials of degree~$1$).

\begin{definition}\label{def:mobius}
Let $M = \left(\begin{smallmatrix}a&b\\c&d\end{smallmatrix}\right)$ be
  a $2 \times 2$ matrix.  Let $f = \sum_{j=0}^{\deg(f)} c_j x^j$ be a
  polynomial in $\ZZ[x]$.  We define the \textbf{M\"obius transform}
  of $f$ induced by $M$ to be the polynomial $\mu_M(f) =
  \sum_{j=0}^{\deg f} c_j \, (ax+b)^j \, (cx+d)^{\deg(f)-j}$.
\end{definition}
In our applications the matrix entries will be integers, and we shall
suppose that at least one of $a$ and $c$ is non-zero.

\begin{definition}\label{def:mobius-degenerate}
A \MT $\mu_M$ is \textbf{degenerate} if $\det M = 0$.
\end{definition}

\begin{definition}\label{def:mobius-inverse}
  Let $\mu_M$ be a M\"obius transform.  We define the
  \textbf{pseudo-inverse} of $\mu_M$ to be the \MT corresponding to
  the classical adjoint $M^{adj} =
  \left(\begin{smallmatrix}d&-b\\-c&a\end{smallmatrix}\right)$.  We
    write $\mu_M^*$ to denote the pseudo-inverse.
\end{definition}

Here is a summary of useful properties of a \MT $\mu_M$.
\begin{proposition}
  \label{prop:moebius}
  Let $M = \left(\begin{smallmatrix}a&b\\c&d\end{smallmatrix}\right)$ be non-singular, so $\mu_M$ is non-degenerate.
\vspace{-6pt}
\begin{enumerate}
  \setlength{\itemsep}{1pt}
  \item[(a)] Let $f = \alpha x+\beta$ be a linear polynomial.  If
    $f(\frac{a}{c}) \neq 0$ 
    then $\mu_M(f)$ is linear; otherwise $\mu_M(f) = \alpha b + \beta
    d$ is a non-zero constant.

  \item[(b)] $\mu_M$ respects multiplication: $\mu_M(gh) = \mu_M(g)\, \mu_M(h)$.

\item [(c)]  $\deg(\mu_M(f)) = \deg(f) \quad \Longleftrightarrow\quad f(\frac{a}{c}) \neq 0$.

\item[(d)] If $\deg(\mu_M(f)) = \deg(f)$ 
  then  $\mu_M^*(\mu_M(f)) = D^{\deg(f)} f(x)$ where $D = \det M$.

\item[(e)] If $\deg(\mu_M^*(f)) = \deg(f)$ 
  then $\mu_M(\mu_M^*(f)) = D^{\deg(f)} f(x)$  where $D = \det M$.

\item[(f)]  If $a,b,c,d \in \ZZ$ and $f \in \ZZ[x]$ is irreducible and $\deg(\mu_M(f)) = \deg(f)$ then $\prim(\mu_M(f))$ is irreducible.
\end{enumerate}
\end{proposition}

\begin{proof}
  Parts~(a) and~(b) are elementary algebra.  Part~(c) follows from~(a) and~(b)
  by considering the factorization of $f$ over a splitting field.
  Parts~(d) and~(e) are elementary for linear $f$; the general case follows
  by repeated application of part~(b).

  For part~(f), suppose we have a counter-example $f \in \ZZ[x]$, then
  we have a non-trivial factorization $\mu_M(f) = gh$, but by~(b) and~(d)
  we deduce that $D^{\deg(f)} f = \mu_M^*(g) \, \mu_M^*(h)$ which is a
  non-trivial factorization, contradicting the assumption that $f$ was
  irreducible.
\end{proof}

Our interest in \MTs is that they offer the possibility of finding a
better LPFW certificate.  Unfortunately we do not yet have a good way
of determining which \MTs are helpful.

\begin{example}
\label{ex:moebius-is-good}
Let $f = 97x^4 +76x^3 +78x^2 +4x +2$.  We obtain a LPFW certificate with
$\rho=7/5$, $\Delta=1$, $n=-4$ with corresponding prime factor $p=10601$.

Let $M = \left(\begin{smallmatrix}1&1\\-3&2\end{smallmatrix}\right)$.
  Let $g = \prim(\mu_M(f)) = (x^4+1)$; by Prop.~\ref{prop:moebius}.(f)
  since $\deg g = \deg f$ a LPFW certificate for $g$ also certifies
  that $f$ is irreducible.  For $g$ we obtain a certificate with
  $\rho=1$, $\Delta=1$, $n=2$ with \textit{much smaller} corresponding
  prime factor $p=17$.
\end{example}

\noindent
\textbf{Unsolved problem: }How to find a good M\"obius matrix $M$ given just $f$?

\subsection{Certifying a transformed polynomial}

Naturally, if we generate a LPFW certificate for a transformed polynomial
$\mu_M(f)$ then we must indicate which \MT was used.
Given two polynomials $f,g \in \ZZ[x]$ of the same degree $d$, and
$M \in \Mat_{2{\times}2}(\ZZ)$, one can easily verify that
$g = \prim(\mu_M(f))$ by evaluating $f$ at $\deg(f)$ distinct rational points, and
$g$ at the (rational) transforms of these points, and then checking
that the ratios of the values are all equal.  So the extra information needed is $M$ and $\mu_M(f)$.

\subsection{Fixed divisors}

\begin{definition}
  Let $f \in \ZZ[x]$ be non-zero.  The \textbf{fixed divisor} of $f$ is
  defined to be $\IntrinsicContent(f) = \gcd \{ f(n) \mid n \in \ZZ \}$.
\end{definition}
Some content-free polynomials have non-trivial fixed divisors:
an example is $f = x^2+x+2$ which is content-free but has fixed divisor~$2$.

\begin{proposition}
  Let $f \in \ZZ[x]$ be non-zero.  Its fixed divisor is equal to:
  \[
  \IntrinsicContent(f) \;=\; \gcd(f(1),f(2),\ldots,f(\deg f))
  \]
\end{proposition}

\begin{proof}
  The standard proof follows easily from representating of $f$ with respect to the
  ``binomial basis'' for $\ZZ[x]$, namely $\{ \binom{x}{k} \mid k \in \NN \}$.
\end{proof}

Polynomials having large fixed divisor $\delta$ cannot have small LPFW
certificates because we are forced to choose large evaluation points
since we must have $(|n|-\rho)^\Delta > \delta$.  This problem becomes
more severe for higher degree polynomials since the fixed divisor can
be as large as $d!$ where $d$ is the degree.


We can reduce the size of the fixed divisor by scaling the
indeterminate (\ie~a \MT for a diagonal matrix), or perhaps
reversing the polynomial and scaling the indeterminate
(\ie~a \MT for an anti-diagonal matrix).  We have not yet
investigated the use of more general \MTs.

Let $f \in \ZZ[x]$ be content-free, irreducible with fixed divisor
$\delta$.  Let $q$ be a prime factor of $\delta$, and let $k$ be the
multiplicity of $q$ in $|f(0)|$.  Then $g(x) = q^{-k} f(q^k x) \in
\ZZ[x]$ has fixed divisor $\delta/q^k$.  In practice, we consider
several polynomials obtained by scaling $x$ by $q^1, q^2,\ldots, q^k$;
in fact scaling by $q^{-1}, q^{-2}, \ldots$ can also be beneficial.

\section{Implementation and Experimentation}

Our prototype implementation runs degree analysis and LPFW search ``in
parallel'': \ie~it repeatedly alternates a few iterations of degree
analysis with a few iterations of LPFW search.  If
degree analysis finds a new factor degree lower bound, $\Delta$, this
information is passed to the LPFW search.

\subsection{Degree analysis}

We adopted the following strategy for choosing primes during degree
analysis: initially we create a list of ``preferential primes''
(\eg~including the first few primes greater than the degree), then we
pick primes alternately from this list or from a random generator.  The range
for randomly generated primes is gradually increased to favour finding
quickly a certificate involving smaller primes (since these are
computationally cheaper to verify).

This strategy was inspired by some experimentation.  There exist
polynomials whose degree analysis certificates must involve ``large''
primes: \eg~a good set of primes for $x^4 +16x^3 +5x^2 -14x -18$ must
contain at least one prime greater than $101$.
Also, empirically we find that a degree analysis certificate for an
(even) Hermite polynomial must use primes greater than the degree.

To issue a certificate, we look for a minimal
cardinality subset of the primes used which suffices.
This subset search is potentially exponential, but in our experiments it is
very rare for a minimal subset to need more than~$3$ primes.

\subsection{Large prime factor witness}

As already mentioned, not all polynomials can be certified irreducible
by degree analysis.  A well-known class of polynomials for which
irreducibility cannot be shown by degree analysis are the
\textit{Swinnerton-Dyer polynomials}: they are the minimal polynomials
for sums of square-roots of ``independent'' integers.  A more general
class of such polynomials was presented
in~\cite{zbMATH03841227}.

We saw in Example~\ref{ex:moebius-is-good}, it can be better to issue
a LPFW certificate for a transformed polynomial, but we do not yet
have a good way of finding a good \MT.  Our current prototype
implementation considers only indeterminate scaling and possibly
reversal: \ie~the M\"obius matrix must be diagonal or anti-diagonal.
A list of all scaling and reverse-scaling transforms by ``simple''
rationals is maintained, and the resulting polynomials are considered
``in parallel''.

For each transformed polynomial we keep track of two evaluation points
(one positive, one negative) and the corresponding evaluations.  The
evaluations are then considered in order of increasing absolute value;
once an evaluation has been processed the corresponding evaluation
point is incremented (or decremented, if it is negative).

The LPFW search depends on a factor degree lower bound, $\Delta$,
which is initially $1$.  The degree analysis ``thread'' may at any
time furnish a better value for $\Delta$.  So that this asynchrony can
work well the LPFW search records, \textit{for each possible factor
  degree lower bound,} any certificates it finds.  When a higher
$\Delta$ is received, the search first checks whether a corresponding
LPFW certificate has already been recorded; if so, that certificate is
produced as output.  Otherwise searching proceeds using the new
$\Delta$.

\subsection{Examples}

Here are a few examples as computed by the current prototype, since
degree analysis picks primes in a pseudo-random order different
certificates may be issued for the same polynomial.
\vspace{-3pt}
\begin{itemize}
  \setlength{\itemsep}{0pt}
\item $x^{16}+4x^{14}+6x^2+4$: degree analysis with prime list $L=[13,127]$ 
\item $x^4 +16x^3 +5x^2 -14x -18$: degree analysis with prime list $L=[107]$ 
\item $21$-st cyclotomic polynomial: LPFW with $\rho=2$, $\Delta=1$, $n=3$, and prime factor $p=368089$  
  \item Swinnerton-Dyer polynomial for $[71,113,163]$: LPFW  with
 $\rho=43$, $\Delta=2$ (with $L=[3]$),  $n=82$ and prime factor $p=2367715751029$  
  \item $97x^4 +76x^3 +78x^2 +4x +2$: \textbf{transform} $x \mapsto \frac{2}{x}$, LPFW $\rho=67/5$, $\Delta=2$ (with $L=[3]$),  $n=-29$ and prime factor $p=3041$ 
\end{itemize}
A quick comment about run-times: our \textit{interpreted prototype} favours
producing certificates which are cheap to verify (rather than cheap to
generate); the degree analysis certificates took $\sim\! 0.25$s each to
generate, the others $\sim\! 0.5$s each.  We did not measure verification
run-time, but fully expect it to be less than 0.01s in each case.
In comparison, the polynomial factorizer in CoCoA took less than 0.01s for
all of these polynomials.

As a larger example: the prototype took $\sim\! 20$s (we expect the
final implementation to be significantly faster) to produce a
certificate for the degree 64 (Swinnerton-Dyer) minimal polynomial of
\[
\sqrt{61}+\sqrt{79}+\sqrt{139}+\sqrt{181}+\sqrt{199}+\sqrt{211}
\]
This polynomial has fixed divisor $\delta = 2^{29} \, 5^{14} \, 13^4
\approx 1.2 \times 10^{28}$.  Our prototype found and applied the transformation
$x \mapsto \frac{52}{15}x$, then produced an LPFW certificate for the transformed
polynomial: $\rho = 451/16$, $\Delta=2$ (with $L=[19]$), $n=46$ and $p
\approx 7.5 \times 10^{180}$ which was confirmed to be ``probably
prime'' (according to GMP~\cite{GMP}).  The classical Berlekamp-Zassenhaus
factorizer in CoCoA~\cite{CoCoA} took about 300s to recognize
irreducibility.

\subsection{A comment about run-time}

An anonymous referee reasonably asked about expected run-time or a
(possibly heuristic) complexity analysis.  The answer is \textit{``It
  depends \ldots''}.  For ``almost all'' polynomials, degree analysis
suffices and is quick.  In our setting, the LPFW search effectively
happens only if a degree analysis certificate cannot be quickly found.
In our experiments, the number of iterations in LPFW search before
producing a certificate was quite irregular.

\section{Conclusion}

As mentioned in the introduction there are many different
criterions for certifying the irreducibility of a polynomial
in $\ZZ[x]$.  Here we have concentrated on just two of them,
and have pointed out how they can ``collaborate''.

We have built a prototype implementation in CoCoA~\cite{CoCoA}, and
plan to integrate it into CoCoALib, the underlying C++ library
(where we expect significant peformance gains).

An interesting future possibility is for the requester of the
certificate to state which criterions may be used (dictated by
the implemented verifiers that the requester has available).
But, a too restrictive choice of criterions may make
it impossible to generate a certificate: \eg~there is no ``Eisenstein''
certificate for most polynomials.


\end{document}